\documentclass[12pt]{amsart}
\usepackage{amsmath,amssymb}

\newtheorem{thm}{Theorem}
\newtheorem{prop}[thm]{Proposition}
\newtheorem{conj}{Conjecture}

\newcommand{\Z}{{\mathbf Z}}
\newcommand{\Q}{{\mathbf Q}}

\newcommand{\C}{{\mathbf C}}
\newcommand{\N}{{\mathbf N}}

\newcommand{\mnorm}{ {\mathcal N}}

\newcommand{\OI}{{\mathcal O}}

\newcommand{\mat}{\operatorname{Mat}}

\begin{document}

\title[Torsion points and common divisors]
{Torsion points on curves and common divisors of $a^k-1$ and $b^k-1$}
\author{ Nir Ailon and Z\'eev Rudnick}
\address{Raymond and Beverly Sackler School of Mathematical Sciences,
Tel Aviv University, Tel Aviv 69978, Israel. Current address:
Department of Computer Science, Princeton University, Princeton, NJ
08544, USA
({\tt nailon@princeton.edu})}
\address{Raymond and Beverly Sackler School of Mathematical Sciences,
Tel Aviv University, Tel Aviv 69978, Israel
({\tt rudnick@post.tau.ac.il})}

\date{December 16, 2002}

\begin{abstract}
We study the behavior of the greatest common divisor of $a^k-1$ and $b^k-1$,
where $a,b$ are fixed integers or polynomials, and $k$ varies. In the
integer case, we conjecture that when $a$ and $b$ are
multiplicatively independent and in addition $a-1$ and $b-1$ are
coprime, then $a^k-1$ and $b^k-1$ are coprime
infinitely often. In the polynomial case, we prove
a strong version of this conjecture. To do this we use a result of
Lang on the finiteness of torsion points on algebraic curves.
We also give a matrix analogue of these results, where for
a nonsingular integer matrix $A$, we look at the greatest common
divisor of the elements of the matrix $A^k-I$.
\end{abstract}

\maketitle


\section{Introduction}

Let $a,b\neq \pm 1$ be nonzero integers.
One of our goals in this paper is to study the common divisors of
$a^k-1$ and $b^k-1$, specifically to understand small values of
$\gcd(a^k-1,b^k-1)$.
If $a=c^u$ and $b=c^v$ for some integer $c$ then clearly $c^k-1$
divides $\gcd(a^k-1,b^k-1)$ and so for the
purpose of understanding small values,
we will assume that $a$ and $b$ are {\em multiplicatively
independent}, that is $a^r\neq b^s$ for $r,s\geq 1$.
Further, since $\gcd(a-1,b-1)$ always divides $\gcd(a^k-1,b^k-1)$, we will
assume that $a-1$ and $b-1$ are coprime.

Based on numerical experiments and other considerations, we conjecture:
\begin{conj}\label{conj-scalar}
If $a,b$ are multiplicatively independent non-zero integers
with $\gcd(a-1,b-1)=1$, then there are infinitely many integers
$k\geq 1$ such that
$$
\gcd(a^k-1,b^k-1) = 1 \;.
$$
\end{conj}

Note that the condition of multiplicative independence of $a$ and
$b$ is not necessary, as the (trivial) example $b=-a$ shows (the
$\gcd$ is $1$ for odd $k$, and $a^k-1$ for even $k$).

A recent result of Bugeaud, Corvaja and Zannier
\cite{BCZ} rules out {\em large} values of $\gcd(a^k-1,b^k-1)$.
They show that if $a,b>1$ are multiplicatively independent positive integers
then for all $\epsilon > 0$,
\begin{equation}\label{bcz}
\gcd(a^k-1,b^k-1) \ll_\epsilon e^{\epsilon k} \; .
\end{equation}
Their argument uses Diophantine approximation techniques and in
particular Schmidt's Subspace Theorem.
They also indicate that there are arbitrarily large values of
$k$ for which the upper bound \eqref{bcz} cannot be significantly improved.

In the function field case, when we replace integers by polynomials,
we are able to prove a strong version of Conjecture \ref{conj-scalar}.

\begin{thm}\label{poly-scal}
Let $f,g\in \C[t]$ be non-constant polynomials.
If $f$ and $g$ are multiplicatively independent, then there exists
a polynomial $h$ such that
\begin{equation}\label{upper bound}
\gcd(f^k-1,g^k-1) \mid h
\end{equation}
for any $k\geq 1$.

If, in addition, $\gcd(f-1,g-1)=1$, then there is a
finite union of proper arithmetic progressions $\cup d_i\N$,
$d_i \geq 2$, such that for $k$
outside these progressions,
$$
\gcd(f^k-1,g^k-1) = 1 \;.
$$
\end{thm}
Note that \eqref{upper bound} is a strong form of \eqref{bcz}.
We derive Theorem \ref{poly-scal} from a result proposed by Lang
\cite{L1} on the finiteness of torsion points on curves - see
section~\ref{sec:Lang}.

We next consider a generalization to the case of matrices.
For an $r\times r$  integer matrix $A\in \mat_r(\Z)$, $A\neq I$,
($I$ being the identity matrix) we define $\gcd(A-I)$
as the greatest common divisor of the entries of $A-I$.
Equivalently, $\gcd(A-I)$ is
the greatest integer $N\geq 1$ such that $A\equiv I \mod N$.
We say that $A$ is {\em primitive} if $\gcd(A-I) = 1$.
Note that $\gcd(A-I)$ divides $\gcd(A^k-I)$ for all $k$.
A similar definition applies to the function field case
$A \in \mat_r(\C[t])$.
We will study behavior of $\gcd(A^k-I)$ as $k$ varies for
a fixed matrix $A$ with coefficients in $\Z$ or in $\C[t]$.
If $\det A=0$ then it holds trivially that $\gcd(A^k-I)=1$ for all
$k\geq 1$. So we will henceforth  assume that $A$ is nonsingular.

For the case of $2\times 2$ matrices, we will show in
section~\ref{sec:2 by 2} that
if $A\in SL_2(\Z)$ is is unimodular and hyperbolic, then $\gcd(A^k-I)$ grows
exponentially as $k\to \infty$.
However, numerical experiments show that for other matrices,
$\gcd(A^k-I)$ displays completely different behaviour.
We formulate the following conjecture:
\begin{conj} \label{conj-at-least-three}
Suppose $r\geq 2$ and $A\in \mat_r(\Z)$ is nonsingular and primitive. Also
assume that there is a pair of eigenvalues of $A$ that are
multiplicatively independent. Then $A^k$ is primitive infinitely
often.
\end{conj}
Note that Conjecture~\ref{conj-at-least-three} subsumes
Conjecture~\ref{conj-scalar}.
It would be interesting to prove an analogue of the upper bound
\eqref{bcz} in this setting.


In section~\ref{section-cyclo} we give an
example where we
can prove Conjecture~\ref{conj-at-least-three}. To describe it,
recall that one may obtain integer matrices by taking an algebraic
integer $u$ in a number field $K$ and letting it act by
multiplication on the ring of integers $\OI_K$ of $K$. This is a
linear map and a choice of integer basis of $\OI_K$ gives us an
integer matrix $A=A(u)$ whose determinant equals the norm of $u$.
We employ
this method for the cyclotomic field $\Q(\zeta_p)$ where $p>3$ is
prime and $\zeta_p$ is a primitive $p$-th root of unity, and $u$
is a non-real unit. We show:
\begin{thm}\label{main-cyc-thm}
Let $u$ be a non-real unit in the extension $\Q(\zeta_p)$, and
$A(u)\in SL_{p-1}(\Z)$ be the corresponding matrix. Then $A(u)^k$
is primitive for all $k\neq 0 \mod p$.
\end{thm}

In the function field case, we have a strong form of
Conjecture~\ref{conj-at-least-three}, which generalizes
Theorem~\ref{poly-scal}:
\begin{thm}\label{poly-mat}
Let $A$ be a nonsingular matrix in $\mat_r(\C[t])$. Assume that either
\begin{enumerate}
\item $A$ is not diagonalizable over the algebraic closure of $\C(t)$,
or
\item $A$ has two eigenvalues that are multiplicatively independent.
\end{enumerate}
Then there exists a polynomial $h$ such that $\gcd(A^k-I)\mid h$ for
any $k$.

If, in addition, $A$ is primitive, then $A^k$ is primitive for all
$k$ outside a finite union of proper arithmetic progressions.
\end{thm}

\emph{Acknowledgements:} We would like to thank Umberto Zannier for
useful discussions and the referee for suggesting several improvements.
Some of the results were part of the first named
author's M.Sc. thesis \cite{Ai} at Tel Aviv University.
The work was partially supported by the Israel Science Foundation,
founded by the Israel Academy of Sciences and Humanities.

\section{Proof of theorem \ref{poly-scal}}\label{sec:Lang}

To prove the theorem, we will use a result which
was conjectured by Serge Lang and proved by Ihara, Serre and Tate (see
\cite{L1} and \cite{L2}), which states that the
intersection of an irreducible curve in $\C^*\times \C^*$ with the roots
of unity $\mu_\infty \times \mu_\infty$ is finite, unless the
curve is of the form $X^nY^m-\zeta=0$ or $X^m-\zeta Y^n=0$ with
$\zeta\in \mu_\infty$, that is unless the curve is the translate of an
algebraic subgroup by a torsion point of $\C^*\times \C^*$. Applying
this result to
the rational curve $\{(f(t), g(t)):t\in \C\}$, we conclude that only
for finitely many
$t$'s both $f(t)$ and $g(t)$ are roots of unity when $f$ and $g$
are multiplicatively independent.

Thus by Lang's theorem we have that there is only a finite set
of points $S \subset \C$ such that for any $s \in S$ both
$f(s)$ and $g(s)$ are roots of unity. So $\gcd(f^k-1,g^k-1)$ can
only have linear factors from $\{(t-s) | s \in S \}$. Write
$$
f^k-1 = \prod_{j=0}^{k-1}(f-\zeta_k^j) \;.
$$
Any two factors on the right side are coprime, so $t-s$ can
divide at most one of them with multiplicity at most $\deg(f)$, and a
similar statement can be said for $g$. Therefore the required
polynomial $h$ can be chosen as
$$ h(t) = \prod_{s \in S} (t-s)^{\min(\deg(f),\deg(g) )} \;.$$

For the second part of theorem \ref{poly-scal}, let $s\in S$ and let
$d_s$ be the least positive integer such that
$$t-s \mid \gcd(f(t)^{d_s}-1, g(t)^{d_s}-1).$$
Then $d_s > 1$ because $\gcd(f-1,g-1)=1$, and clearly for $k\notin
d_s\N$,
$$t-s \nmid \gcd(f(t)^k-1,g(t)^k-1).$$
Then $\cup_{ s \in S} d_s \N $ is the required finite union of proper
arithmetic progressions outside of which $\gcd(f^k-1,g^k-1)=1$. \qed

Note that Theorem~\ref{poly-mat} implies Theorem~\ref{poly-scal}. We
have chosen to give the proof of Theorem~\ref{poly-scal} separately to
illustrate the ideas in a simple context.

\section{$2\times 2$ matrices}\label{sec:2 by 2}
Let $A\in SL_2(\Z)$ be a $2\times 2$ unimodular
matrix which is {\em hyperbolic}, that is $A$ has two distinct
real eigenvalues. We show:
\begin{prop}\label{two-by-two}
Let $A\in SL_2(\Z)$ be a hyperbolic matrix with eigenvalues $\epsilon,
\epsilon^{-1}$, where $|\epsilon|>1$.
Then $\gcd(A^k-I) \gg |\epsilon|^{k/2}$.
\end{prop}
\begin{proof}
\footnote{We thank the referee for suggesting this proof, which
replaces our original, more complicated, version.}
Let $K$ be the real quadratic field $\Q(\epsilon)$ and $\OI_K$ its ring
of integers.   We may diagonalize the matrix $A$ over $K$, that is
write $A=P\begin{pmatrix}\epsilon&0\\0&\epsilon^{-1} \end{pmatrix}P^{-1}$
with $P$ a $2\times 2$ matrix having entries in $K$. Since $P$ is only
determined up to a scalar multiple, we may, after multiplying $P$ by an
algebraic integer of $\OI_K$, assume that $P$ has entries in $\OI_K$.
Then $P^{-1} =\frac 1{\det(P)} P^{ad}$
where $P^{ad}$ also has entries  in $\OI_K$. Thus
we have
$$
A^k-I = \frac 1{\det(P)}
P\begin{pmatrix}\epsilon^k-1&0\\0&\epsilon^{-k}-1 \end{pmatrix}P^{ad} \;.
$$

The entries of $ A^k-I $ are thus $\OI_K$-linear combinations
$(\epsilon^k-1)/\det(P)$ and of $(\epsilon^{-k}-1)/\det(P)$.
We now note that
$$\epsilon^{-k}-1 = -\epsilon^{-k}(\epsilon^k-1)$$
and thus the entries of $ A^k-I $ are all $\OI_K$-multiples of
$(\epsilon^k-1)/\det(P)$. In particular, $\gcd(A^k-I)$, which is a $\Z$-linear
combination of the entries of $A^k-I$, can be written as
$$
\gcd(A^k-I) = \frac {\epsilon^k-1}{\det(P)} \gamma_k
$$
with $\gamma_k\in \OI_K$.

Now taking norms from $K$ to $\Q$ we see
$$
|\gcd(A^k-I) |^2  =
\frac{|\mnorm(\epsilon^k-1)| }{|\mnorm(\det P)|}|\mnorm(\gamma_k)|
\;.
$$
Since $\gamma_k\neq 0$, we have $|\mnorm(\gamma_k)|\geq 1$ and thus
$$
|\gcd(A^k-I) |^2 \geq \frac{ |\mnorm(\epsilon^k-1)|} {|\mnorm(\det
 P)|}
\gg \epsilon^k
$$
which gives $|\gcd(A^k-I) | \gg \epsilon^{k/2}$.
\end{proof}

A special case of this Proposition appeared as a problem in the 54-th
W.L. Putnam Mathematical Competition, 1994, see  \cite[pages 82,
242]{Andreescu}.

\section{Cyclotomic Fields}\label{section-cyclo}

A standard construction of unimodular matrices is to take a unit
$u$ of norm one in a number field $K$ and let it act by
multiplication on the ring of integers $\OI_K$ of $K$. This gives
a linear map and a choice of integer basis of $\OI_K$ gives us an
integer matrix whose determinant equals the norm of $u$ and is
thus unimodular. We employ this method for the case when $u$ is a
nonreal unit to give a construction of matrices $A$ with the
property that $A^k$ is primitive infinitely often.

We recall some basic facts on units in a cyclotomic field.  Let
$p>3$ be a prime, $\zeta_p$ a primitive $p$-th root of unity, and
$K=\Q(\zeta_p)$ the cyclotomic extension of the rationals.  It is
a field of degree $p-1$.  The ring of integers of this field
$\OI_K$ is $\Z[\zeta_p]$.  $K$ is purely imaginary, therefore the
norm function is positive, and the norm of a unit $u$ is always
$1$. Also note that the structure of the unit group $E_p$ of
$\OI_K$ is:
\begin{equation}\label{group struct}
E_p = W_pE_p^+,
\end{equation}
 where $W_p$ are the roots of
unity in $K$ and $E_p^+$ is the group of the real units in
$\OI_K$. A proof of this fact can be found, for example, in
\cite[Theorem 4.1]{L3}.

\subsection{Proof of Theorem~\ref{main-cyc-thm}}
We now prove theorem~\ref{main-cyc-thm}, that is show that if
$u\in E_p\backslash E_p^+$ is a non-real unit and $k\not\equiv
0\mod p$ then the matrix corresponding to $u^k$ is primitive.

The method we will use is that if we choose a basis
$\omega_0=1,\omega_1,\dots,\omega_{p-2}$ of $\Z[\zeta_p]$ and take a unit
$U$ in $\Z[\zeta_p]$, then we get a matrix $A(U)=(a_{i,j})$ whose
entries are determined by
$$
U\omega_i=\sum_{j=0}^{p-2} a_{j,i} \omega_j \;.
$$
In particular if we find that in the expansion of
$$U=U\cdot \omega_0 = \sum_{j=0}^{p-2} a_{j,0} \omega_j
$$
we have an index $j\neq 0$ so that $a_{j,0}=a_{0,0}$, then in the
matrix $A(U)-I$ corresponding to $U-1$, the first column will contain the
entries $a_{0,0}-1$ and $a_{j,0}=a_{0,0}$ which are clearly coprime
and thus the matrix $A(U)$ is {\em primitive}.

Another option is to have $a_{0,0}=0$ in which case in the matrix
of $U-1$, the $(0,0)$ entry is $-1$ and thus again $A(U)$ is
primitive. We will apply this method to the case that $U=u^k$ is a
power of a non-real unit $u$ and $k\not\equiv 0\mod p$.


Let $u\in E_p\backslash E_p^+$ is a non-real unit. By \eqref{group
struct}, we can write:
$$ u = \zeta_p^x u^+ $$
where $u^+$ is a {\em real} unit and $x$ is an integer not
congruent to $0$ mod $p$. Therefore,
$$
u^k = \zeta_p^{xk} (u^+)^k
$$
and
$$
\zeta_p^{-xk} u^k = (u^+)^k
$$
is real. Therefore it can be represented as an integer combination of
$\zeta_p, \zeta_p^2, \dots, \zeta_p^{p-1}$ as follows:
$$\zeta_p^{-xk} u^k = \sum_{j=1}^{p-1}\alpha_j \zeta_p^j $$
where $\alpha_j = \alpha_{p-j}$ for each $j$. For convenience we will
set $\alpha_0:=0$.

Multiplying by $\zeta_p^{xk}$, we find
$$ u^k = \sum_{j=0}^{p-1}\alpha_j \zeta_p^{j+xk} $$
and changing the summation variable,
$$ u^k = \sum_{i=0}^{p-1}\alpha_{i-xk} \zeta_p^i $$
where the index of $\alpha$ is calculated mod $p$.
Using the relation
$$ \zeta_p^{p-1} = -1 - \zeta_p -\dots - \zeta_p^{p-2} $$
we find that in terms of the integer basis $\omega_j=\zeta_p^{j}$,
$j=0,\dots,p-2$ we have
$$
u^k = \sum_{i=0}^{p-2}(\alpha_{i-xk} - \alpha_{p-1-xk})\omega_i \;.
$$

If $k\not\equiv 0\mod p$ then $2xk\not\equiv 0\mod p$ since
$x\not\equiv 0\mod p$. If $2xk \not\equiv -1\mod p$ then the
coefficients of $\omega_0$ and $\omega_{2xk}$ are
equal. Therefore $u^k$ is primitive. If $2xk$ is congruent to $-1$
mod $p$, then
the coefficient of $\omega_0$ vanishes and thus in this case as well,
$u^k$ is primitive.

Thus we found that if $k\not\equiv 0 \mod p$, the matrix corresponding
to $u^k$ is primitive.

\qed

Note that by virtue of \eqref{group struct}, the eigenvalues of
$A(u)$ come in complex conjugate pairs whose ratios are $p$-th
roots of unity. This is somewhat similar to the trivial scalar example
described in the introduction, namely $b=\pm a$.


\section{Proof of theorem \ref{poly-mat}}
We extend the idea of the proof of Theorem~\ref{poly-scal}
to cover the matrix case.
We first show that there is only a finite set $S$ of points $s\in
\C$ such that $t-s$ divides $\gcd(A^k-I)$ for some $k$.

Let $M$ be a matrix such that $MAM^{-1}$ is in
Jordan form. The elements of $M$ are meromorphic functions on the
Riemann surface $R$ corresponding to some finite extension of
$\C(t)$. Denote by $pr:R\to {\mathbb  P}^1$ the associated projection of $R$ to the projective line.
Let $S_0$ be the finite collection of poles of these
functions.

Assume first that $A$ is not diagonalizable over the algebraic
closure of $C(t)$.
 Thus for any $t_0 \in R \backslash S_0$, $A(t_0)$ is
not diagonalizable, and therefore $A(t_0)^k-I \neq
0$ for all $k$ (recall that a matrix of finite order ($A^m=I$)
is automatically diagonalizable), in other words, $(t-t_0)$ does not divide
$\gcd(A^k-I)$. Thus only the finitely many linear forms $t-s$,
where $s\in pr(S_0)$ is the projection of some point in $S_0$, can divide
$\gcd(A^k-I)$.

We denote by $\lambda_i(t)$ the eigenvalues of $A$ which are
multivalued functions of $t$, that is meromorphic functions on the
Riemann surface.  Assume now that  $\lambda_1$ and $\lambda_2$ are
multiplicatively independent, and that $A$ is diagonalizable.
Suppose that $(t-t_0) \mid \gcd(A^k-I)$ for some $k>1$ and $t_0
\in R \backslash S_0$. Then $A^k-I$ evaluated at $t_0$ is the zero
matrix, and also:
$$
M(t_0)(A(t_0)^k-I)M(t_0)^{-1} = 0 \;,
$$
and we deduce that
$$
\lambda_1(t_0)^k-1 = \lambda_2(t_0)^k-1 = 0\;.
$$
In particular, $\lambda_1(t_0)$ and $\lambda_2(t_0)$ are roots of
unity. Thus, we reduce to proving that $\lambda_1$ and $\lambda_2$
can be simultaneous roots of unity only at a finite set of points.

To prove this, we want to use Lang's theorem for the curve in
$\C^2$ parameterized by $(\lambda_1(t), \lambda_2(t))$.
%
Denote by $Y$ the Zariski closure of the image of the map
$(\lambda_1,\lambda_2):R\backslash S_0\to \C^2$.
$Y$ is an irreducible algebraic curve in $\C^2$.
If $Y$ is of dimension $0$, then it is a point, so $\lambda_1(t)$ and
$\lambda_2(t)$ are constants,
and since they are multiplicatively independent none of them
can be a root of unity. Otherwise, we may apply Lang's
theorem for this curve and conclude that unless the curve $Y$ is
of the form $F^m-\zeta G^n=0$ or $F^mG^n=\zeta$ with $\zeta$ a
root of unity (which is not the case when $\lambda_1$ and
$\lambda_2$ are multiplicatively independent) , it has only
finitely many torsion points. In other words, there can only be
finitely many points of the form $(\zeta_1, \zeta_2)$ on $Y$,
where $\zeta_1$ and $\zeta_2$ are roots of unity.

We now prove that there is a polynomial $h$ such that $\gcd(A^k-I)$
divides $h$ for all $k$. Since there is a finite set $S$ of possible zeros
of $\gcd(A^k-I)$, it suffices to show that the multiplicity of a zero
of $\gcd(A^k-I)$ is bounded.

Write $B =MAM^{-1}$, so $B$ is in Jordan form. Denote by $v_{t_0}(f)$ the
multiplicity of the zero at $t_0 \in R$ of $f$. So clearly, for
any $t_0 \in R$ there exists $c(t_0) \in \N$ such that
$$ v_{t_0}(\gcd(A^k-I)) \leq c(t_0) + v_{t_0}(\gcd(B^k-I)), $$
and for all $t_0$ outside the finite set $S_0$ of poles of entries
of $M$, $c(t_0) = 0$. So it suffices to prove that
$v_{t_0}(\gcd(B^k-I))$ is bounded.

Clearly,
$$
\gcd(B^k-I) \mid \det(B^k-I) =
\prod_{j=0}^{k-1}\det(B-\zeta_k^jI) \;,
$$
where $\zeta_k$ is a primitive $k$-th root of unity.
Denoting the diagonal elements of
$B-I$ by $b_1,...,b_r$, 
we see that
$$
\det(B^k-I) = \prod_{d=1}^r\prod_{j=0}^{k-1} (b_d - \zeta_k^j) \;.
$$

Because a meromorphic function on a Riemann surface has a finite
degree,
reasoning as in the proof of theorem \ref{poly-scal} we see that
for any $t_0 \in R$, $v_{t_0}(\prod_{j=1}^k(b_d-\zeta_k^j))$ is
bounded, for all $k$. Therefore $v_{t_0}(\det(B^k-I))$ is bounded
for all $k$.

Now assume in addition that $A$ is primitive: $\gcd(A-I) = 1$.
For any $s\in S$, the set of $k$'s such
that $A(s)^k=I$, i.e. $(t-s)\mid \gcd(A^k-I)$,
is an arithmetic progression $d_s\Z$ which is
proper since it does not contain $1$.
Therefore the set of $k$ such that $\gcd(A^k-I) \neq 1$ is a finite union
of proper arithmetic progressions, and hence for $k$ outside this finite
union of proper arithmetic progressions, we have $\gcd(A^k-I) = 1$. \qed


\begin{thebibliography}{99}
\bibitem[Ai]{Ai}
Ailon, N. Primitive powers of matrices and related problems, Tel Aviv
University M.Sc Thesis, October 2001.

\bibitem[An]{Andreescu}
Andreescu, T, and  Gelca, R.
{\em Mathematical Olympiad challenges}.
Birkhauser Boston, Inc., Boston, MA, 2000.

\bibitem[BCZ]{BCZ}
Bugeaud, Y., Corvaja, P. and Zannier, U. An upper bound for the G.C.D
of $a^n-1$ and $b^n-1$, to appear in Math. Zeitschrift.


\bibitem[L1]{L1}
Lang, S. Annali di Matematica pura ed applicata (IV), Vol. LXX,
1965 229--234.

\bibitem[L2]{L2}
Lang, S. Fundamentals of Diophantine Geometry. Springer-Verlag
1983 200--207.

\bibitem[L3]{L3}
Lang, S. Cyclotomic Fields, Springer-Verlag 1978 79--82.


\bibitem[W]{W}
Washington, L. C. Introduction to Cyclotomic Fields.
Springer-Verlag 1982 29--38,143--146.

\end{thebibliography}
\end{document}